# An Extension to Gaussian Semigroup and Some Applications [*]

Guibao Liu [#]

*Abstract.* We look at the semigroup generated by a system of heat equations. Applications to testing normality and option pricing are addressed.

*Keywords.* Heat equation, Semigroup, Normality test, Option pricing

## 1 Extension to Gaussian Semigroup

It is well-known that Gaussian semigroup solves the heat equation $u_t = u_{xx}$. Consider breaking the heat equation into pieces, each with its own domain and diffusivity.

**Proposition 1.** *Solution to the system of partial differential equations*

(1.1) $\quad \dfrac{\partial u_k}{\partial t} = \dfrac{1}{2}\sigma_k^2 \dfrac{\partial^2 u_k}{\partial x^2} \qquad q_k < x < q_{k-1};\ k = 1,2,\ldots,N;\ q_0 = +\infty,\ q_N = -\infty$

(1.2) $\quad u_l = f_l(x) = \delta(x) \qquad$ at $t = 0$ (initial condition, if phase $l$ brackets $x = 0$; $\delta(x)$ is the Dirac-delta function)

(1.3) $\quad u_{k \ne l} = f_k(x) = 0 \qquad$ at $t = 0$ (initial condition)

(1.4) $\quad u_k = u_{k+1} \qquad$ at $x = q_k$ (boundary condition, $k = 1,2,\ldots,N\text{-}1$)

(1.5) $\quad \dfrac{1}{2}\sigma_k^2 \dfrac{\partial u_k}{\partial x} = \dfrac{1}{2}\sigma_{k+1}^2 \dfrac{\partial u_{k+1}}{\partial x} \qquad$ at $x = q_k$ (boundary condition, $k = 1,2,\ldots,N\text{-}1$)

*where* $u_k = u_k(x,t)$, *is unique, represents a probability distribution and is a semigroup.*

*Proof.* (i) (Uniqueness) Suppose $w$ is another solution. Then $v := w - u$ solves the system with initial conditions that are identically zero. For $F := \int_{-\infty}^{+\infty} \dfrac{1}{2} v^2 dx$,

$\dfrac{dF}{dt} = \int_{-\infty}^{+\infty} v v_t dx = \sum_{k=1}^{N} \int_{q_k}^{q_{k-1}} \dfrac{1}{2}\sigma_k^2 v_k v_{kxx} dx = \sum_{k=1}^{N} \dfrac{1}{2}\sigma_k^2 v_{kx} v_k \Big|_{q_k}^{q_{k-1}} - \sum_{k=1}^{N} \dfrac{1}{2}\sigma_k^2 \int_{q_k}^{q_{k-1}} \left(\dfrac{\partial v_k}{\partial x}\right)^2 dx$. Invoke the regularity conditions that $v$ and $v_x$ vanish at $\pm\infty$; the first summation is zero by (1.4) and (1.5), so $\dfrac{dF}{dt} \le 0$. At $t = 0$, $F = 0$. By definition, $F$ is non-negative. So, $F = 0\ \forall t$, i.e. $v = 0\ \forall t$.

---

(ii) (Probability distribution) By the maximum principle, $u_k$ is non-negative for all $k$. Let
$$F = \int_{-\infty}^{+\infty} u\,dx \text{; then } \frac{dF}{dt} = \int_{-\infty}^{+\infty} \frac{\partial u}{\partial t}\,dx = \sum_{k=1}^{N} \int_{q_k}^{q_{k-1}} \frac{1}{2}\sigma_k^2 u_{kxx}\,dx = \sum_{k=1}^{N} \frac{1}{2}\sigma_k^2 u_{kx}\Big|_{q_k}^{q_{k-1}} = 0.$$
Since $F(x,0) = \int_{-\infty}^{+\infty} \delta(x)\,dx = 1$, $F = 1\ \forall t$.

(iii) (Semigroup) We show that $u$ satisfies the Chapman-Kolmogorov Equation,

(1.6) $$u(x,t;\{q_i\}) = \int_{-\infty}^{\infty} u(y,s;\{q_i\})u(x-y,t-s;\{q_i-y\})\,dy$$

where $\{\cdot\}$ denotes the set of interphase boundaries. For $q_k < x < q_{k-1}$, (1.6) implies

(1.7) $$u_k(x,t;\{q_i\}) = \int_{-\infty}^{+\infty} u(y,s;\{q_i\})u_k(x-y,t-s;\{q_i-y\})\,dy \qquad (k=1,2,...,N).$$

Denote the convolution in (1.7) as $w_k(x,t)$; we have

(1.8) 
$$\frac{\partial}{\partial t} w_k = \int_{-\infty}^{\infty} u(y,s;\{q_i\}) \frac{\partial}{\partial t} u_k(x-y,t-s;\{q_i-y\})\,dy$$
$$= \frac{1}{2}\sigma_k^2 \int_{-\infty}^{\infty} u(y,s;\{q_i\}) \frac{\partial^2}{\partial x^2} u_k(x-y,t-s;\{q_i-y\})\,dy = \frac{1}{2}\sigma_k^2 \frac{\partial^2}{\partial x^2} w_k$$

(1.9)
$$w_k(q_k,t) = \int_{-\infty}^{\infty} u(y,s;\{q_i\})u_k(q_k-y,t-s;\{q_i-y\})\,dy$$
$$= \int_{-\infty}^{\infty} u(y,s;\{q_i\})u_{k+1}(q_k-y,t-s;\{q_i-y\})\,dy = w_{k+1}(q_k,t)$$

(1.10)
$$\frac{1}{2}\sigma_k^2 \frac{\partial}{\partial x} w_k(q_k,t) = \int_{-\infty}^{\infty} u(y,s;\{q_i\}) \frac{1}{2}\sigma_k^2 \frac{\partial}{\partial x} u_k(q_k-y,t-s;\{q_i-y\})\,dy$$
$$= \int_{-\infty}^{\infty} u(y,s;\{q_i\}) \frac{1}{2}\sigma_{k+1}^2 \frac{\partial}{\partial x} u_{k+1}(q_k-y,t-s;\{q_i-y\})\,dy.$$
$$= \frac{1}{2}\sigma_{k+1}^2 \frac{\partial}{\partial x} w_{k+1}(q_k,t)$$

(1.11) $$w_k(x,0) = \int_{-\infty}^{\infty} \delta(y)\cdot 0\,dy = 0 \qquad \text{if } k \neq l \qquad \text{and}$$

(1.12) $$w_k(x,0) = \int_{-\infty}^{\infty} \delta(y)\delta(x-y)\,dy = \delta(x) \qquad \text{if } k = l.$$

Hence $w_k(x,t)$, $k = 1,2,...,N$ solve the system (1.1)-(1.5). By uniqueness, $w(x,t)$ equals $u(x,t)$. (1.6) shows that $u$ is a convolution semigroup which reduces to Gaussian if $\sigma_k = \sigma\ \forall k$.

**Proposition 2 (Two-phase Distribution).** *Solution to the system (1.1)-(1.5) when $N = 2$ ( let $q_1 = q$ ) is*

$$(1.13) \quad u_1(x,t) = \frac{2\sigma_1}{\sigma_1 + \sigma_2} \frac{1}{\sqrt{2\pi}\sigma_1 \sqrt{t}} \exp\left\{-\frac{1}{2}\left[\frac{x - \left(1 - \frac{\sigma_1}{\sigma_2}\right)q}{\sigma_1 \sqrt{t}}\right]^2\right\} \quad (q < x < \infty)$$

$$(1.14) \quad u_2(x,t) = \frac{1}{\sqrt{2\pi}\sigma_2 \sqrt{t}} \exp\left[-\frac{1}{2}\left(\frac{x}{\sigma_2 \sqrt{t}}\right)^2\right] + \frac{\sigma_2 - \sigma_1}{\sigma_1 + \sigma_2} \frac{1}{\sqrt{2\pi}\sigma_2 \sqrt{t}} \exp\left[-\frac{1}{2}\left(\frac{x - 2q}{\sigma_2 \sqrt{t}}\right)^2\right]$$

$$(-\infty < x < q)$$

*if $q > 0$, and is*

$$(1.15) \quad u_1(x,t) = \frac{1}{\sqrt{2\pi}\sigma_1 \sqrt{t}} \exp\left[-\frac{1}{2}\left(\frac{x}{\sigma_1 \sqrt{t}}\right)^2\right] + \frac{\sigma_1 - \sigma_2}{\sigma_1 + \sigma_2} \frac{1}{\sqrt{2\pi}\sigma_1 \sqrt{t}} \exp\left[-\frac{1}{2}\left(\frac{x - 2q}{\sigma_1 \sqrt{t}}\right)^2\right]$$

$$(q < x < \infty)$$

$$(1.16) \quad u_2(x,t) = \frac{2\sigma_2}{\sigma_1 + \sigma_2} \frac{1}{\sqrt{2\pi}\sigma_2 \sqrt{t}} \exp\left\{-\frac{1}{2}\left[\frac{x - \left(1 - \frac{\sigma_2}{\sigma_1}\right)q}{\sigma_2 \sqrt{t}}\right]^2\right\} \quad (-\infty < x < q)$$

*if $q < 0$.* [1]

*Proof.* See Appendix A.

**Proposition 3 (Three-phase Distribution).** *Solution to the system (1.1)-(1.5) when $N = 3$ and $l = 2$ is*

$$(1.17) \quad u_1(x,t) = \frac{\sigma_1}{\sigma_2} \cdot \frac{1}{\sqrt{2\pi}\sigma_1 \sqrt{t}} \sum_{n=-\infty}^{+\infty}\left\{\exp\left[-\frac{1}{2}\left(\frac{x - q_1 + (\sigma_1/\sigma_2)|(2n+1)q_1 - 2nq_2|}{\sigma_1 \sqrt{t}}\right)^2\right]\right.$$

$$\left. + \exp\left[-\frac{1}{2}\left(\frac{x - q_1 + (\sigma_1/\sigma_2)|(2n-1)q_1 - 2nq_2|}{\sigma_1 \sqrt{t}}\right)^2\right]\right\}$$

$$(0 < q_1 < x < \infty)$$

$$(1.18) \quad u_2(x,t) = u_{2,1} - u_{2,2} - u_{2,3} \quad (q_2 < x < q_1)$$

---

[1] Border case like $q = 0$ ($x = q$) can be assigned to either solution.

where

$$u_{2,1}(x,t) = \frac{1}{\sqrt{2\pi}\sigma_2\sqrt{t}} \sum_{n=-\infty}^{+\infty} \left\{ \exp\left[-\frac{1}{2}\left(\frac{x+2nq_1-2nq_2}{\sigma_2\sqrt{t}}\right)^2\right] + \exp\left[-\frac{1}{2}\left(\frac{x+2nq_1-(2n+2)q_2}{\sigma_2\sqrt{t}}\right)^2\right] \right\}$$

$$u_{2,2}(x,t) = \frac{\sigma_3}{\sigma_2} \cdot \frac{1}{\sqrt{2\pi}\sigma_2\sqrt{t}} \sum_{n=-\infty}^{+\infty}\sum_{m=-\infty}^{+\infty} \left\{ \exp\left[-\frac{1}{2}\left(\frac{|x+2nq_1-(2n+1)q_2|+|2mq_1-(2m-1)q_2|}{\sigma_2\sqrt{t}}\right)^2\right] \right.$$

$$\left. + \exp\left[-\frac{1}{2}\left(\frac{|x+2nq_1-(2n+1)q_2|+|2mq_1-(2m+1)q_2|}{\sigma_2\sqrt{t}}\right)^2\right] \right\}$$

$$u_{2,3}(x,t) = \frac{\sigma_1}{\sigma_2} \cdot \frac{1}{\sqrt{2\pi}\sigma_2\sqrt{t}} \sum_{n=-\infty}^{+\infty}\sum_{m=-\infty}^{+\infty} \left\{ \exp\left[-\frac{1}{2}\left(\frac{|x+(2n-1)q_1-2nq_2|+|(2m+1)q_1-2mq_2|}{\sigma_2\sqrt{t}}\right)^2\right] \right.$$

$$\left. + \exp\left[-\frac{1}{2}\left(\frac{|x+(2n-1)q_1-2nq_2|+|(2m-1)q_1-2mq_2|}{\sigma_2\sqrt{t}}\right)^2\right] \right\},$$

and

(1.19)
$$u_3(x,t) = \frac{\sigma_3}{\sigma_2} \cdot \frac{1}{\sqrt{2\pi}\sigma_3\sqrt{t}} \sum_{n=-\infty}^{+\infty} \left\{ \exp\left[-\frac{1}{2}\left(\frac{x-q_2-(\sigma_3/\sigma_2)|2nq_1-(2n-1)q_2|}{\sigma_3\sqrt{t}}\right)^2\right] \right.$$

$$\left. + \exp\left[-\frac{1}{2}\left(\frac{x-q_2-(\sigma_3/\sigma_2)|2nq_1-(2n+1)q_2|}{\sigma_3\sqrt{t}}\right)^2\right] \right\}$$

$$(-\infty < x < q_2 < 0).$$

*Proof.* See Appendix B.

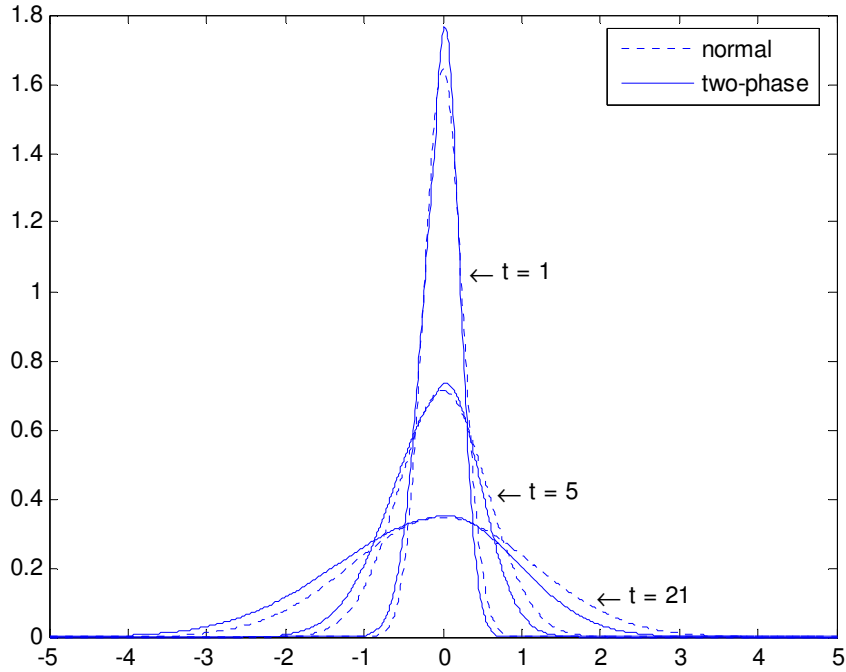

Figure 1: Two-phase versus Normal

Note: $\sigma_1 = 0.2$, $\sigma_2 = 0.3$, $q = -0.1$.

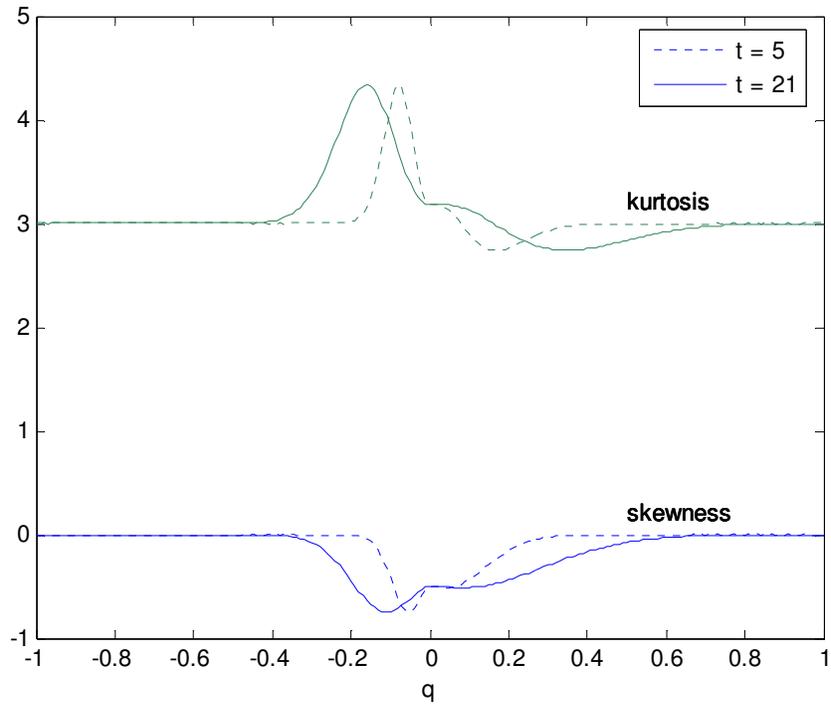

Figure 2: Skewness and kurtosis of the two-phase distribution

Note: $\sigma_1 = 2.5\%$, $\sigma_2 = 5\%$. Dash-lined: $t = 5$; Solid-lined: $t = 21$.

Figure 1 compares the two-phase distribution with the normal distribution of commensurate variance for various scales in $t$. Figure 2 shows variation in skewness and kurtosis across different values of $q$ given both $\sigma_1$ and $\sigma_2$.

## 2 Some Applications

*2.1 Testing Normality*

The two-phase distribution derived above offers a simple way to test normal distribution based on maximum-likelihood inference. Using the null hypothesis of $\sigma_1 = \sigma_2$, the likelihood ratio test statistic follows $\chi^2$ distribution with one degree of freedom. Log likelihood function to the two-phase distribution looks as follows:

(2.1)
$$\log L = n\left(\log 2 - \log(\sigma_1 + \sigma_2) - \frac{1}{2}\log 2\pi\right) - \frac{1}{2}\sum_{(i)=1}^{n}\left(\frac{x_{(i)} - \left(1 - \frac{\sigma_1}{\sigma_2}\right)q}{\sigma_1 \sqrt{t}}\right)^2$$

$$+ \sum_{(i)=n+1}^{N} \log\left(\frac{1}{\sqrt{2\pi}\sigma_2 \sqrt{t}} e^{-\frac{1}{2}\left(\frac{x_{(i)}}{\sigma_2 \sqrt{t}}\right)^2} + \frac{\sigma_2 - \sigma_1}{\sigma_1 + \sigma_2}\frac{1}{\sqrt{2\pi}\sigma_2 \sqrt{t}} e^{-\frac{1}{2}\left(\frac{x_{(i)} - 2q}{\sigma_2 \sqrt{t}}\right)^2}\right)$$

if $q > 0$, and

(2.2)
$$\log L = \sum_{(i)=1}^{n} \log\left(\frac{1}{\sqrt{2\pi}\sigma_1 \sqrt{t}} e^{-\frac{1}{2}\left(\frac{x_{(i)}}{\sigma_1 \sqrt{t}}\right)^2} + \frac{\sigma_1 - \sigma_2}{\sigma_1 + \sigma_2}\frac{1}{\sqrt{2\pi}\sigma_1 \sqrt{t}} e^{-\frac{1}{2}\left(\frac{x_{(i)} - 2q}{\sigma_1 \sqrt{t}}\right)^2}\right)$$

$$+ (N - n)\left(\log 2 - \log(\sigma_1 + \sigma_2) - \frac{1}{2}\log 2\pi\right) - \frac{1}{2}\sum_{(i)=n+1}^{N}\left(\frac{x_{(i)} - \left(1 - \frac{\sigma_2}{\sigma_1}\right)q}{\sigma_2 \sqrt{t}}\right)^2$$

if $q \leq 0$.

In (2.1)-(2.2), $N$ is sample size, $x_{(i)}$ represents the $i^{th}$ largest observation by value, and $n$ denotes the total number of observations that are greater than $q$ (i.e., $x_{(N)} \leq ... \leq x_{(n+1)} \leq q \leq x_{(n)} \leq ... \leq x_{(1)}$). Nonlinearity requires estimating parameters $\sigma_1, \sigma_2$ and $q$ with a numerical procedure. Standard Newton's method usually suffices the purpose. Table 1 reports the ML parameter estimates and test statistics for a variety of samples of weekly returns on the S&P 500 index (SPX).

Table 1: Maximum likelihood estimates and likelihood ratio tests: two-phase vs. normal

| sample | N | $\hat{\sigma}_1$ | $\hat{\sigma}_2$ | $\hat{q}$ | $\hat{\sigma}$ | LR | $p$-value |
|---|---|---|---|---|---|---|---|
| 1/2/1996 – 12/27/2005 | 521 | 0.990 | 3.552 | -2.377 | 1.038 | 26.461 | 2.69E-7 |
|  |  | (0.031) | (1.233) | (0.108) | (0.032) |  |  |
| 1/5/2004 – 12/27/2005 | 104 | 0.594 | 0.817 | -0.981 | 0.611 | 0.723 | 0.395 |
|  |  | (0.044) | (0.315) | (0.190) | 0.042 |  |  |
| 1/7/2002 – 12/29/2003 | 104 | 1.023 | 1.890 | -1.807 | 1.084 | 3.066 | 0.080 |
|  |  | (0.075) | (0.735) | (0.235) | (0.075) |  |  |
| 1/3/2000 – 12/31/2001 | 104 | 1.204 | 3.961 | -2.298 | 1.357 | 11.258 | 7.93E-4 |
|  |  | (0.086) | (1.801) | (0.229) | (0.094) |  |  |
| 1/5/1998 – 12/27/1999 | 104 | 1.170 | 1.054 | -0.001 | 1.116 | 0.843 | 0.359 |
|  |  | (0.103) | (0.096) | (0.120) | (0.077) |  |  |
| 1/2/1996 – 12/29/1997 | 105 | 0.971 | 0.400 | -1.032 | 0.870 | 6.981 | 8.24E-3 |
|  |  | (0.094) | (0.142) | (0.192) | (0.060) |  |  |

Note: (1) $N$: sample size; $\hat{\sigma}$: ML estimate to the standard deviation of normal distribution; $LR$: likelihood ratio test statistic. (2) Parameter estimates are in percentage. (3) Standard errors are in parentheses. (4) Convergence criterion is set at 10E-8.

## 2.2 Pricing European Option

Following standard option pricing theory, the two-phase distribution can be used to define a risk neutral probability measure to price derivative securities like European options. For equities, assume that time-$t$ stock price $S_t$ has the dynamics

$$(2.3) \quad S_t = S_0 \exp\left(\int_0^t \mu_s ds + Z_t\right)$$

where $Z_t$ is a stochastic process implied by the two-phase distribution.[2] The drift term $\int_0^t \mu_s ds$ can be resolved by requiring $S_t$ normalized by money market account as a martingale. For simplicity, we use a constant risk-free interest rate ($r$). Applying transformation of distribution, it is easy to obtain a risk-neutral measure $Q$ so that $S_0 = e^{-rt} E^Q(S_t)$. With this measure, a European call option $c$ with strike $K$ is priced straightforwardly by evaluating $e^{-r(T-t)} E^Q\left((S_T - K)^+ | S_t\right)$. We present the solution in Equations (2.4)-(2.7).

$$(2.4) \quad c = S\Psi_{11} - Ke^{-r(T-t)}\Psi_{12} \quad \text{if } 0 < q < -\ln(S/K) - \bar{\mu}(T-t),$$
$$(2.5) \quad c = S\Psi_{21} - Ke^{-r(T-t)}\Psi_{22} \quad \text{if } q > 0 \vee \left(-\ln(S/K) - \bar{\mu}(T-t)\right)$$
$$(2.6) \quad c = S\Psi_{31} - Ke^{-r(T-t)}\Psi_{32} \quad \text{if } q < 0 \wedge \left(-\ln(S/K) - \bar{\mu}(T-t)\right)$$
$$(2.7) \quad c = S\Psi_{41} - Ke^{-r(T-t)}\Psi_{42} \quad \text{if } -\ln(S/K) - \bar{\mu}(T-t) < q < 0$$

---

[2] $Z_0 = 0$.

where

$$\Psi_{11} = \frac{2\sigma_1}{\sigma_1 + \sigma_2} e^{(\bar{\mu} + \frac{1}{2}\sigma_1^2 - r)(T-t) + (1 - \frac{\sigma_1}{\sigma_2})q} \Phi\left(\frac{\ln(S/K) + (\bar{\mu} + \sigma_1^2)(T-t) + (1 - \frac{\sigma_1}{\sigma_2})q}{\sigma_1\sqrt{T-t}}\right)$$

$$\Psi_{12} = \Phi\left(\frac{\ln(S/K) + \bar{\mu}(T-t) + (1 - \frac{\sigma_1}{\sigma_2})q}{\sigma_1\sqrt{T-t}}\right)$$

$$\Psi_{21} = e^{(\bar{\mu} + \frac{1}{2}\sigma_2^2 - r)(T-t)} \left[\Phi\left(\frac{\ln(S/K) + (\bar{\mu} + \sigma_2^2)(T-t)}{\sigma_2\sqrt{T-t}}\right) - \Phi\left(\frac{-q + \sigma_2^2(T-t)}{\sigma_2\sqrt{T-t}}\right) + \frac{\sigma_2 - \sigma_1}{\sigma_1 + \sigma_2} e^{2q} \cdot \left(\Phi\left(\frac{\ln(S/K) + (\bar{\mu} + \sigma_2^2)(T-t) + 2q}{\sigma_2\sqrt{T-t}}\right) - \Phi\left(\frac{q + \sigma_2^2(T-t)}{\sigma_2\sqrt{T-t}}\right)\right)\right]$$
$$+ e^{(\bar{\mu} + \frac{1}{2}\sigma_1^2 - r)(T-t) + (1 - \frac{\sigma_1}{\sigma_2})q} \Phi\left(\frac{-q + \sigma_1\sigma_2(T-t)}{\sigma_2\sqrt{T-t}}\right)$$

$$\Psi_{22} = \frac{2\sigma_2}{\sigma_1 + \sigma_2} \Phi\left(\frac{\ln(S/K) + \bar{\mu}(T-t)}{\sigma_2\sqrt{T-t}}\right) + \frac{\sigma_1 - \sigma_2}{\sigma_1 + \sigma_2} \Phi\left(-\frac{\ln(S/K) + \bar{\mu}(T-t) + 2q}{\sigma_2\sqrt{T-t}}\right)$$

$$\Psi_{31} = e^{(\bar{\mu} + \frac{1}{2}\sigma_1^2 - r)(T-t)} \Phi\left(\frac{\ln(S/K) + (\bar{\mu} + \sigma_1^2)(T-t)}{\sigma_1\sqrt{T-t}}\right)$$
$$+ \frac{\sigma_1 - \sigma_2}{\sigma_1 + \sigma_2} e^{(\bar{\mu} + \frac{1}{2}\sigma_1^2 - r)(T-t) + 2q} \Phi\left(\frac{\ln(S/K) + (\bar{\mu} + \sigma_1^2)(T-t) + 2q}{\sigma_1\sqrt{T-t}}\right)$$

$$\Psi_{32} = \Phi\left(\frac{\ln(S/K) + \bar{\mu}(T-t)}{\sigma_1\sqrt{T-t}}\right) + \frac{\sigma_1 - \sigma_2}{\sigma_1 + \sigma_2} \Phi\left(\frac{\ln(S/K) + \bar{\mu}(T-t) + 2q}{\sigma_1\sqrt{T-t}}\right)$$

$$\Psi_{41} = \frac{2\sigma_2}{\sigma_1 + \sigma_2} e^{(\bar{\mu} + \frac{1}{2}\sigma_2^2 - r)(T-t) + (1 - \frac{\sigma_2}{\sigma_1})q} \left[\Phi\left(\frac{\ln(S/K) + (\bar{\mu} + \sigma_2^2)(T-t) + (1 - \frac{\sigma_2}{\sigma_1})q}{\sigma_2\sqrt{T-t}}\right) - \Phi\left(\frac{-q + \sigma_1\sigma_2(T-t)}{\sigma_1\sqrt{T-t}}\right)\right]$$
$$+ e^{(\bar{\mu} + \frac{1}{2}\sigma_1^2 - r)(T-t)} \Phi\left(\frac{-q + \sigma_1^2(T-t)}{\sigma_1\sqrt{T-t}}\right) + \frac{\sigma_1 - \sigma_2}{\sigma_1 + \sigma_2} e^{(\bar{\mu} + \frac{1}{2}\sigma_1^2 - r)(T-t) + 2q} \Phi\left(\frac{q + \sigma_1^2(T-t)}{\sigma_1\sqrt{T-t}}\right)$$

$$\Psi_{42} = \frac{2\sigma_2}{\sigma_1+\sigma_2}\Phi\left(\frac{\ln(S/K)+\bar{\mu}(T-t)+(1-\frac{\sigma_2}{\sigma_1})q}{\sigma_2\sqrt{T-t}}\right)+\frac{\sigma_1-\sigma_2}{\sigma_1+\sigma_2}$$

$\Phi(\cdot)$ is the standard normal cdf, $\bar{\mu} = \mu(T-t)/(T-t)$, $\mu(t) = rt - \ln\Lambda(t)$

and

$$\Lambda(t) = e^{\frac{1}{2}\sigma_1^2 t}\left[\frac{2\sigma_1}{\sigma_1+\sigma_2}e^{(1-\frac{\sigma_1}{\sigma_2})q}\Phi\left(\frac{-q+\sigma_1\sigma_2 t}{\sigma_2\sqrt{t}}\right)\right]+e^{\frac{1}{2}\sigma_2^2 t}\left[\Phi\left(\frac{q-\sigma_2^2 t}{\sigma_2\sqrt{t}}\right)+\frac{\sigma_2-\sigma_1}{\sigma_1+\sigma_2}e^{2q}\Phi\left(\frac{-q-\sigma_2^2 t}{\sigma_2\sqrt{t}}\right)\right]$$

if $q > 0$ and

$$\Lambda(t) = e^{\frac{1}{2}\sigma_1^2 t}\left[\Phi\left(\frac{-q+\sigma_1^2 t}{\sigma_1\sqrt{t}}\right)+\frac{\sigma_1-\sigma_2}{\sigma_1+\sigma_2}e^{2q}\Phi\left(\frac{q+\sigma_1^2 t}{\sigma_1\sqrt{t}}\right)\right]+e^{\frac{1}{2}\sigma_2^2 t}\left[\frac{2\sigma_2}{\sigma_1+\sigma_2}e^{(1-\frac{\sigma_2}{\sigma_1})q}\Phi\left(\frac{q-\sigma_1\sigma_2 t}{\sigma_1\sqrt{t}}\right)\right]$$

if $q < 0$.

Table 2 compares the two-phase estimates with the classical Black-Scholes model prices. Figure 3 and Figure 4 depict implied volatility curves and surface, respectively, by the two-phase formula.

Table 2: A comparison of call prices based on two-phase-distribution (1) versus Black-Scholes call prices (2)

| | \multicolumn{8}{c|}{$S=100$, $r=5\%$ annual} |
|---|---|---|---|---|---|---|---|---|
| | (1) | (2) | (1) | (2) | (1) | (2) | (1) | (2) |
| $\tau$ | $K=80$ | | $K=85$ | | $K=90$ | | $K=95$ | |
| 17  | 20.192 | 20.189 | 15.252 | 15.234 | 10.507 | 10.461 | 6.304  | 6.263  |
| 45  | 20.673 | 20.627 | 16.046 | 15.969 | 11.801 | 11.715 | 8.128  | 8.082  |
| 80  | 21.474 | 21.381 | 17.166 | 17.056 | 13.262 | 13.168 | 9.86   | 9.824  |
| 136 | 22.838 | 22.708 | 18.861 | 18.737 | 15.258 | 15.174 | 12.074 | 12.062 |
| 227 | 24.95  | 24.809 | 21.294 | 21.185 | 17.962 | 17.911 | 14.97  | 15     |
| 318 | 26.882 | 26.755 | 23.434 | 23.353 | 20.271 | 20.259 | 17.4   | 17.473 |
| $\tau$ | $K=100$ | | $K=105$ | | $K=110$ | | $K=115$ | |
| 17  | 3.094  | 3.139  | 1.157  | 1.282  | 0.319  | 0.422  | 0.065 | 0.112 |
| 45  | 5.173  | 5.217  | 3.005  | 3.146  | 1.586  | 1.772  | 0.761 | 0.935 |
| 80  | 7.023  | 7.079  | 4.775  | 4.93   | 3.096  | 3.321  | 1.918 | 2.169 |
| 136 | 9.335  | 9.417  | 7.045  | 7.226  | 5.191  | 5.455  | 3.739 | 4.057 |
| 227 | 12.324 | 12.45  | 10.023 | 10.248 | 8.055  | 8.372  | 6.402 | 6.791 |
| 318 | 14.821 | 14.991 | 12.53  | 12.8   | 10.516 | 10.882 | 8.767 | 9.214 |

Notes: $S$ is spot price of stock, $K$ is call strike, $r$ is risk-free rate, and $\tau$ is time-to-expiration measured in days. The call prices are based on the parameterization $[\sigma_1, \sigma_2, q] = [0.3, 0.4, -0.02]$.

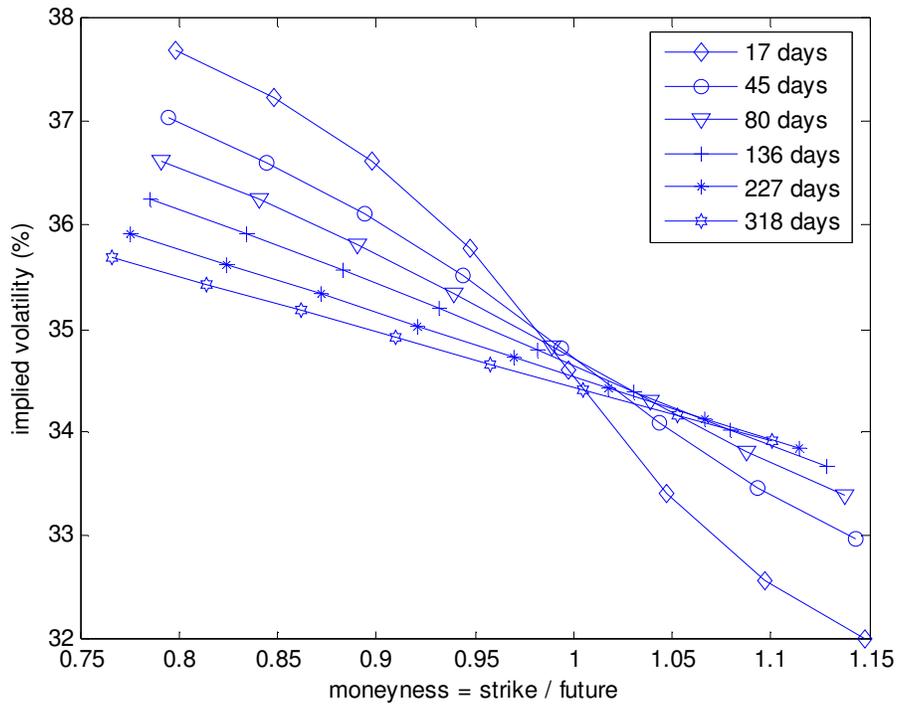

Figure 3: Implied volatilities using the two-phase model with $[\sigma_1, \sigma_2, q] = [0.3, 0.4, -0.02]$

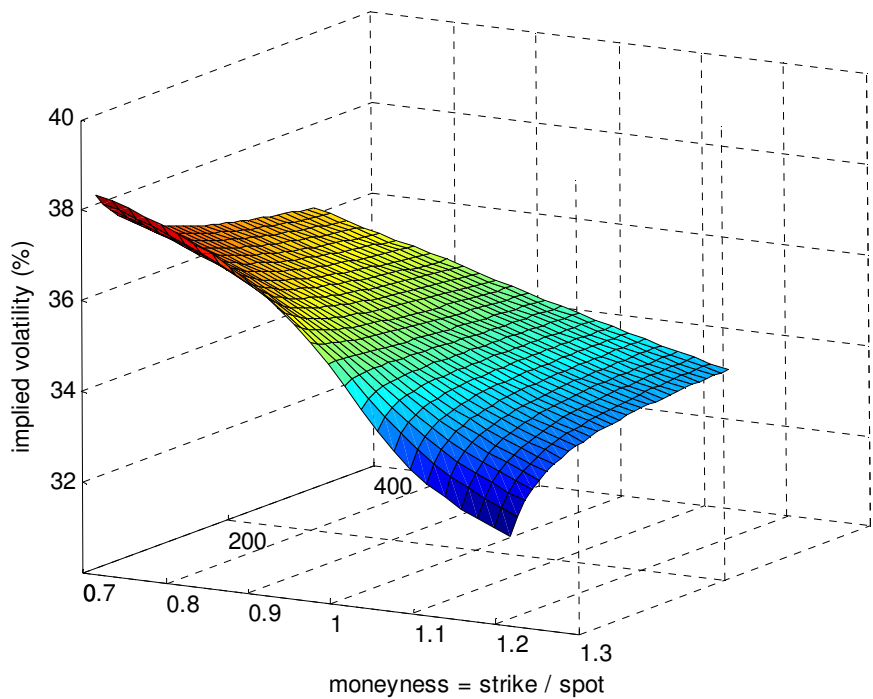

Figure 4: Implied volatility surface using the two-phase model with $[\sigma_1, \sigma_2, q] = [0.3, 0.4, -0.02]$

# Appendix A

Without loss of generality, let $q > 0$. Through change of variables $y = x - q$ and $v(y,t) = w(x,t)$, (1.1)-(1.5) become:

(A1) $$\frac{\partial v_1}{\partial t} = a_1 \frac{\partial^2 v_1}{\partial y^2} \qquad 0 < y < \infty$$

(A2) $$\frac{\partial v_2}{\partial t} = a_2 \frac{\partial^2 v_2}{\partial y^2} \qquad -\infty < y < 0$$

(A3) $$v_1(y,0) = 0$$

(A4) $$v_2(y,0) = \delta(y+q)$$

(A5) $$v_1(0,t) = v_2(0,t)$$

(A6) $$a_1 v_{1y}(0,t) = a_2 v_{2y}(0,t)$$

where $a_i = \frac{1}{2}\sigma_i^2$ $i = 1,2$. A semi-analytical solution for the system (A1)-(A6) is given in Polyanin (2002): [3]

---

[3] See also Luikov (1968).

(A7) $$v_1(y,t) = -\frac{1}{\sqrt{\pi a_1}} \int_0^t \exp\left(-\frac{y^2}{4a_1(t-\tau)}\right) \frac{g(\tau)}{\sqrt{t-\tau}} d\tau$$

(A8) $$v_2(y,t) = \frac{1}{2\sqrt{\pi a_2 t}} \left( \exp\left(-\frac{(y+q)^2}{4a_2 t}\right) + \exp\left(-\frac{(y-q)^2}{4a_2 t}\right) \right)$$
$$+ \frac{1}{\sqrt{\pi a_2}} \int_0^t \exp\left(-\frac{y^2}{4a_2(t-\tau)}\right) \frac{g(\tau)}{\sqrt{t-\tau}} d\tau$$

where

(A9) $$g(t) = -\frac{\sqrt{a_1}}{\pi(\sqrt{a_1}+\sqrt{a_2})} \frac{d}{dt} \int_0^t \exp\left(-\frac{q^2}{4a_2\tau}\right) \frac{1}{\sqrt{\tau}\sqrt{t-\tau}} d\tau.$$

It can be shown that

(A10) $$\int_0^t \exp\left(-\frac{q^2}{4a_2\tau}\right) \frac{1}{\sqrt{\tau}\sqrt{t-\tau}} d\tau = \pi \cdot erfc\left(\frac{|q|}{2\sqrt{a_2 t}}\right)$$

where $erfc(\cdot)$ is the complementary error function.

When $a_1 = a_2$, (A1)-(A6) reduce to the heat equation $\frac{\partial v}{\partial t} = a_2 \frac{\partial^2 v}{\partial y^2}$ ($-\infty < y < \infty$) with solution

(A11) $$v(y,t) = \frac{1}{2\sqrt{\pi a_2 t}} \exp\left(-\frac{(y+q)^2}{4a_2 t}\right).$$

Equate (A7) and (A8) to (A11); we find,

(A12) $$\int_0^t \exp\left(-\frac{y^2}{4a_2(t-\tau)}\right) \exp\left(-\frac{q^2}{4a_2\tau}\right) \frac{1}{\tau\sqrt{\tau}\sqrt{t-\tau}} d\tau = \frac{2\sqrt{\pi a_2}}{|q|\sqrt{t}} \exp\left(-\frac{(|y|+|q|)^2}{4a_2 t}\right).$$

Define $A := \int_0^t \exp\left(-\frac{y^2}{4a_2(t-\tau)}\right) \exp\left(-\frac{q^2}{4a_2\tau}\right) \frac{1}{\tau\sqrt{\tau}\sqrt{t-\tau}} d\tau$ and

$B := \int_0^t \exp\left(-\frac{y^2}{4a_2(t-\tau)}\right) \exp\left(-\frac{q^2}{4a_2\tau}\right) \frac{1}{\sqrt{\tau}\sqrt{t-\tau}} d\tau$ and use the fact that $A = \frac{-2a_2}{q} \frac{dB}{dq}$,

(A13) $$B = \pi \cdot erfc\left(\frac{|y|}{2\sqrt{a_2 t}} + \frac{|q|}{2\sqrt{a_2 t}}\right).$$

This suggests a symmetry solution to the following integral:

(A14) $$\int_0^t \exp\left(-\frac{\alpha^2}{(t-\tau)}\right)\exp\left(-\frac{\beta^2}{\tau}\right)\frac{1}{\sqrt{\tau}\sqrt{t-\tau}}d\tau = \pi \cdot erfc\left(\frac{|\alpha|}{\sqrt{t}} + \frac{|\beta|}{\sqrt{t}}\right).$$

Let $C := \int_0^t \exp\left(-\frac{y^2}{4a_1(t-\tau)}\right)\exp\left(-\frac{q^2}{4a_2\tau}\right)\frac{1}{\sqrt{\tau}\sqrt{t-\tau}}d\tau$ and

$D := \int_0^t \exp\left(-\frac{y^2}{4a_1(t-\tau)}\right)\exp\left(-\frac{q^2}{4a_2\tau}\right)\frac{1}{\tau\sqrt{\tau}\sqrt{t-\tau}}d\tau$; by the identity $D = \frac{-2a_2}{q}\frac{dC}{dq}$,

(A15) $$D = \frac{2\sqrt{\pi a_2}}{|q|\sqrt{t}}\exp\left(-\left(\frac{|y|}{2\sqrt{a_1 t}} + \frac{|q|}{2\sqrt{a_2 t}}\right)^2\right).$$

Substituting (A12) and (A15) to (A7) and (A8), it is straightforward to obtain the solution of (1.13)-(1.14). By symmetry, Equations (1.15)-(1.16) follow when $q < 0$.

## Appendix B

We solve for (1.17)-(1.19) using Green's functions. In general, solution to heat equation $v_t = av_{xx}$ in the domain of $0 \le x \le l$ and with Neumann boundary conditions can be expressed as

(B4) $$v(x,t) = \int_0^l f(\xi)G_1(x,\xi,t)d\xi - a\int_0^t g_2(\tau)G_1(x,0,t-\tau)d\tau + a\int_0^t g_1(\tau)G_1(x,l,t-\tau)d\tau$$

where $f(x)$ represents the initial condition, $g_1(t) = \frac{\partial}{\partial x}v(l,t)$, $g_2(t) = \frac{\partial}{\partial x}v(0,t)$ and the Green's function $G_1(x,\xi,t)$ is:

(B5) $$G_1(x,\xi,t) = \frac{1}{2\sqrt{\pi at}}\sum_{n=-\infty}^{+\infty}\left\{\exp\left[-\frac{(x-\xi+2nl)^2}{4at}\right] + \exp\left[-\frac{(x+\xi+2nl)^2}{4at}\right]\right\}.$$

Likewise, for $v_t = av_{xx}$ in the domain of $0 \le x < \infty$ or $-\infty < x \le 0$, the solutions are,

(B6) $$v(x,t) = \int_0^\infty f(\xi)G_2(x,\xi,t)d\xi - a\int_0^t g_1(\tau)G_2(x,0,t-\tau)d\tau$$

and

(B7) $$v(x,t) = \int_{-\infty}^{0} f(\xi)G_2(x,\xi,t)d\xi + a\int_0^t g_1(\tau)G_2(x,0,t-\tau)d\tau$$

respectively, where $f(x) = v(x,0)$, $g_1(t) = \dfrac{\partial}{\partial x}v(0,t)$ and $G_2$ has the expression,

(B8) $$G_2(x,\xi,t) = \frac{1}{2\sqrt{\pi at}}\left\{\exp\left[-\frac{(x-\xi)^2}{4at}\right] + \exp\left[-\frac{(x+\xi)^2}{4at}\right]\right\}.\ [4]$$

In reference to (1.1)-(1.5), let $\dfrac{1}{2}\sigma_1^2 u_{1x}(q_1,t) = \dfrac{1}{2}\sigma_2^2 u_{2x}(q_1,t) = g_1(t)$ and $\dfrac{1}{2}\sigma_2^2 u_{2x}(q_2,t) = \dfrac{1}{2}\sigma_3^2 u_{3x}(q_2,t) = g_2(t)$. With proper changes in variable and domain, it can be shown that

(B9) $$u_1(x,t) = \int_0^\infty f_1(\xi+q_1)\frac{1}{\sqrt{2\pi}\sigma_1\sqrt{t}}\left\{\exp\left[-\frac{(x-\xi-q_1)^2}{2\sigma_1^2 t}\right] + \exp\left[-\frac{(x+\xi-q_1)^2}{2\sigma_1^2 t}\right]\right\}d\xi$$
$$-\int_0^t g_1(\tau)\frac{2}{\sqrt{2\pi}\sigma_1\sqrt{t-\tau}}\exp\left(-\frac{(x-q_1)^2}{2\sigma_1^2(t-\tau)}\right)d\tau$$
$$(x > q_1)$$

(B10) $$u_2(x,t) = \int_0^{q_1-q_2} f_2(\xi+q_2)\frac{1}{\sqrt{2\pi}\sigma_2\sqrt{t}}\sum_{n=-\infty}^{+\infty}\left\{\exp\left[-\frac{(x-\xi+2nq_1-(2n+1)q_2)^2}{2\sigma_2^2 t}\right]\right. +$$
$$\left.\exp\left[-\frac{(x+\xi+2nq_1-(2n+1)q_2)^2}{2\sigma_2^2 t}\right]\right\}d\xi - \int_0^t g_2(\tau)\frac{2}{\sqrt{2\pi}\sigma_2\sqrt{t-\tau}} \cdot$$
$$\sum_{n=-\infty}^{+\infty}\exp\left[-\frac{(x+2nq_1-(2n+1)q_2)^2}{2\sigma_2^2(t-\tau)}\right]d\tau + \int_0^t g_1(\tau)\frac{2}{\sqrt{2\pi}\sigma_2\sqrt{t-\tau}} \cdot$$
$$\sum_{n=-\infty}^{+\infty}\exp\left[-\frac{(x+(2n-1)q_1-2nq_2)^2}{2\sigma_2^2(t-\tau)}\right]d\tau$$
$$(q_2 < x < q_1)$$

and

(B11) $$u_3(x,t) = \int_{-\infty}^{0} f_3(\xi+q_2)\frac{1}{\sqrt{2\pi}\sigma_3\sqrt{t}}\left\{\exp\left[-\frac{(x-\xi-q_2)^2}{2\sigma_3^2 t}\right] + \exp\left[-\frac{(x+\xi-q_2)^2}{2\sigma_3^2 t}\right]\right\}d\xi$$
$$+\int_0^t g_2(\tau)\frac{2}{\sqrt{2\pi}\sigma_3\sqrt{t-\tau}}\exp\left(-\frac{(x-q_2)^2}{2\sigma_3^2(t-\tau)}\right)d\tau$$

---

[4] See Polyanin (2002) or Kevorkian (1990).

$$(x < q_2),$$

where $f_i$ represents the initial condition in the $i$ th phase. Evaluating (B9)-(B11) at the interphase boundaries, we have:

(B12) $$u_1(q_1,t) = -\frac{2}{\sqrt{2\pi}\sigma_1} \int_0^t \frac{g_1(\tau)}{\sqrt{t-\tau}} d\tau$$

(B13)
$$u_2(q_1,t) = \frac{1}{\sqrt{2\pi}\sigma_2\sqrt{t}} \sum_{n=-\infty}^{+\infty} \left\{ \exp\left[-\frac{((2n+1)q_1 - 2nq_2)^2}{2\sigma_2^2 t}\right] + \exp\left[-\frac{((2n-1)q_1 - 2nq_2)^2}{2\sigma_2^2 t}\right] \right\}$$
$$-\int_0^t g_2(\tau) \frac{2}{\sqrt{2\pi}\sigma_2 \sqrt{t-\tau}} \sum_{n=-\infty}^{+\infty} \exp\left(-\frac{((2n+1)q_1 - (2n+1)q_2)^2}{2\sigma_2^2(t-\tau)}\right) d\tau$$
$$+\int_0^t g_1(\tau) \frac{2}{\sqrt{2\pi}\sigma_2 \sqrt{t-\tau}} \sum_{n=-\infty}^{+\infty} \exp\left(-\frac{(2nq_1 - 2nq_2)^2}{2\sigma_2^2(t-\tau)}\right) d\tau$$

(B14)
$$u_2(q_2,t) = \frac{1}{\sqrt{2\pi}\sigma_2\sqrt{t}} \sum_{n=-\infty}^{+\infty} \left\{ \exp\left[-\frac{(2nq_1 - (2n-1)q_2)^2}{2\sigma_2^2 t}\right] + \exp\left[-\frac{(2nq_1 - (2n+1)q_2)^2}{2\sigma_2^2 t}\right] \right\}$$
$$-\int_0^t g_2(\tau) \frac{2}{\sqrt{2\pi}\sigma_2 \sqrt{t-\tau}} \sum_{n=-\infty}^{+\infty} \exp\left(-\frac{(2nq_1 - 2nq_2)^2}{2\sigma_2^2(t-\tau)}\right) d\tau$$
$$+\int_0^t g_1(\tau) \frac{2}{\sqrt{2\pi}\sigma_2 \sqrt{t-\tau}} \sum_{n=-\infty}^{+\infty} \exp\left(-\frac{((2n+1)q_1 - (2n+1)q_2)^2}{2\sigma_2^2(t-\tau)}\right) d\tau$$

and

(B15) $$u_3(q_2,t) = \frac{2}{\sqrt{2\pi}\sigma_3} \int_0^t \frac{g_2(\tau)}{\sqrt{t-\tau}} d\tau.$$

By continuity, $u_1(q_1,t) = u_2(q_1,t)$ and $u_2(q_2,t) = u_3(q_2,t)$. We find two Abel's integral equations of the form

(B16) $$\int_0^t \frac{h(\tau)}{\sqrt{t-\tau}} d\tau = j(t).$$

Solution to Equation (B16) is,

(B17) $$h(t) = \frac{1}{\pi} \frac{\partial}{\partial t} \int_0^t \frac{j(\tau)}{\sqrt{t-\tau}} d\tau.$$

Since

$$\lim_{\tau \to t} \sum_{n=-\infty}^{+\infty} \exp\left(-\frac{((2n+1)q_1 - (2n+1)q_2)^2}{2\sigma_2^2(t-\tau)}\right) = \lim_{\tau \to t} \lim_{N \to \infty} \sum_{n=-N}^{N} \exp\left(-\frac{((2n+1)q_1 - (2n+1)q_2)^2}{2\sigma_2^2(t-\tau)}\right),$$

$$= \lim_{N \to \infty} \lim_{\tau \to t} \sum_{n=-N}^{N} \exp\left(-\frac{((2n+1)q_1 - (2n+1)q_2)^2}{2\sigma_2^2(t-\tau)}\right) = 0$$

and by using the identity (A10), we obtain

(B18)
$$g_1(t) = -\frac{\sigma_1}{2\sqrt{\pi}\sigma_2} \sum_{n=-\infty}^{+\infty} \left[ \exp\left(-\frac{((2n+1)q_1 - 2nq_2)^2}{2\sigma_2^2 t}\right) \frac{|(2n+1)q_1 - 2nq_2|}{\sqrt{2}\sigma_2 t^{3/2}} \right.$$
$$\left. + \exp\left(-\frac{((2n-1)q_1 - 2nq_2)^2}{2\sigma_2^2 t}\right) \frac{|(2n-1)q_1 - 2nq_2|}{\sqrt{2}\sigma_2 t^{3/2}} \right].$$

Similarly,

(B19)
$$g_2(t) = \frac{\sigma_3}{2\sqrt{\pi}\sigma_2} \sum_{n=-\infty}^{+\infty} \left[ \exp\left(-\frac{(2nq_1 - (2n-1)q_2)^2}{2\sigma_2^2 t}\right) \frac{|2nq_1 - (2n-1)q_2|}{\sqrt{2}\sigma_2 t^{3/2}} \right.$$
$$\left. + \exp\left(-\frac{(2nq_1 - (2n+1)q_2)^2}{2\sigma_2^2 t}\right) \frac{|2nq_1 - (2n+1)q_2|}{\sqrt{2}\sigma_2 t^{3/2}} \right].$$

Finally, substituting (B19) and (B20) in (B9)-(B11) and applying (A14) to simplify, Equations (B1)-(B3) are obtained after some algebra.

**Disclaimer**